\documentclass[leqno,a4paper,11pt]{article}
\usepackage{amsmath,amssymb,amsthm}

\newtheorem{theorem}{Theorem}[section]
\newtheorem{lem}[theorem]{Lemma}
\newtheorem{prop}[theorem]{Proposition}
\newtheorem{cor}[theorem]{Corollary}

\title{On filling minimality of simple Finsler manifolds.}
\author{\tt henrik.koehler@ruhr-uni-bochum.de}

\begin{document}

\maketitle

\begin{abstract}
\noindent
This paper states a formula for the difference of the Holmes-Thompson volumes of two simple Finsler manifolds of arbitrary dimension, in terms of the difference of the boundary distances and their derivatives.
An application is a conditioned result on filling minimality.
\end{abstract}

\section{Introduction.}

Let $M$ be a smooth compact manifold with boundary $\partial M$ and a reversible Finsler metric $F$.
$(M,F)$ is called \emph{simple}, if it is convex, without conjugate points, and any two points $x,y\in M$ are connected by a unique geodesic segment.
Simple manifolds are known to be contractible, and whether a manifold is simple can be determined from the data of boundary distances (see \cite{Cr2}; the transfer from Riemannian to Finsler metrics has no influence).

In this article, $(M,F)$ shall be called \emph{minimal (Finsler volume) filling}, if $\mathrm{vol}_{\tilde F}(\tilde M)\ge\mathrm{vol}_F(M)$ holds for all oriented Finsler manifolds $(\tilde M,\tilde F)$ with $\partial\tilde M=\partial M$ and $\mathrm{dist}_{\tilde F}(y,z)\ge\mathrm{dist}_F(y,z)$ $\forall\,y,z\in\partial M$, where $\mathrm{vol}$ denotes the Holmes-Thompson (sc.~symplectic) volume.
The notion of filling volume was originally introduced in \cite{Gr} in the context of systolic and isoperimetric inequalities.
It should be mentioned, that the Holmes-Thompson volume coincides with the standard volume in the Riemannian case; hence the above notion comprises filling minimality for Riemannian manifolds.

An open question is, whether simple manifolds are minimal fillings.
In contrast, a manifold that contains regions which are not (or too sparsely) intersected by minimal geodesics between boundary points, clearly cannot be a minimal filling.
Therefore, some restriction has to be imposed on $(M,F)$ to guarantee that the data of boundary distances give sufficient information about the interior of $M$; here simplicity seems a capable requirement.

In the Riemannian case, the question of filling minimality is often considered together with the boundary rigidity question, which asks, whether a Riemannian metric is determined (up to isometries) from its boundary distances.
Filling minimality was proved for conformal metrics and for two-dimensional Riemannian SGM-manifolds (see \cite{CrDa}), and for metrics close to one another in a $C^{3,\alpha}$-sense (see \cite{CrDaSh}).
In two recent articles (\cite{BuIv1} resp.~\cite{BuIv2}), the problems of boundary rigidity and filling minimality were solved for simple Riemannian metrics close to the flat resp.~hyperbolic metric in a $C^2$ resp.~$C^3$-sense.
Also, filling minimality was recently shown for two dimensional Finsler metrics with minimal geodesics (see \cite{Iv1}).
Further, a local result was obtained in \cite{Iv2}, stating volume monotonicity w.r.t.~boundary distance increasing changes of the Finsler metric in a $C^\infty$-neighbourhood for simple Finsler manifolds of any dimension.

This article states in cor.~\ref{minfilcor}, that an inequality for the boundary distances of two simple Finsler manifolds implies the same inequality between the symplectic volumes,\smallskip\\
\begin{tabular}{cl}
if&the dimension is $n=2$ (as already known from \cite{Iv1}),\\
or&$n=3$ or $n=4$ and the sum of the boundary distances is again\\
&a boundary distance function of some simple Finsler manifold,\\
or&the boundary distance functions are $C^2$-close to each other.\\%,\\
%or&there is a connecting $C^1$-homotopy of simple Finsler metrics\\
%&with increasing boundary distances.
\end{tabular}\par
It should be noticed, that the third condition needs no assumption (other than simplicity) on the metrics in the interior; thereby it differs from results like prop.~1.2 in \cite{CrDaSh} or thm.~2 of \cite{Iv2} on volume monotonicity w.r.t.~small changes of the Riemann resp.~Finsler metric.
To clarify what ``$C^2$-close'' means for boundary distances, their behaviour near the diagonal is examined in section \ref{diag.s}.
One might ask, whether the second condition is necessary; however, prop.~\ref{ctex} shows, that for $n=3$, the sum of boundary distance functions need not come from a simple Finsler manifold.

The essential tool is a relationship between the canonical symplectic two-form on the co-tangent bundle and boundary distances (cf.~\cite{Ot}).
This allows to represent the boundary integral in Santal\'o's formula in terms of the mixed second derivative of the boundary distance function (see prop.~\ref{santalo2}).
Using this identity, prop.~\ref{minfil} expresses the difference of Finsler volumes as an integral of the difference of boundary distances; thereby it generalizes what was known for two-dimensional Riemannian manifolds (thm.~1.4 of \cite{CrDa}).

\textsc{Acknowledgements:}
The author is grateful to Sergei Ivanov for his comments on a prior preprint version.
He also kindly provided a proof for $C^{1,1}$-regularity of the exponential map along the zero section in the Finsler case (prop.~\ref{C11}).
Further, the author would like to thank Christopher Croke and Gerhard Knieper for their helpful remarks.

\section{Santal\'o-type integral formulas.}

In all what follows, only simple Finsler manifolds are considered.
Since these are always contractible, one may restrict to the model case of an $n$-disk.

Henceforth, let $B=\{x\in\mathbb R^n:\|x\|<1\}$ denote the unit ball, $S^{n-1}$ its boundary and $\bar B=B\cup S^{n-1}$ its closure.
Suppose $\bar B$ is equipped with a reversible Finsler metric $F:T\bar B\to[0,\infty)$, i.e.~$F$ is a norm on every $T_x\bar B$, depending smoothly on $x\in\bar B$, $F(-v)=F(v)\;\forall\,v\in T\bar B$, and the bilinear form associated to $F$ at $w\in T_x\bar B\setminus\{0\}$ via
$$g_w(u,v):=\frac{d^2}{2ds\,dt}\bigg|_{s=t=0}F^2(w+su+tv)
\quad(u,v\in T_x\bar B),$$
is positive definite on $T_x\bar B$.
For later use, notice that $g_w(w,w)=F^2(w)$ and $g_{rw}=g_w\;\forall\,r\not=0$.
Further, let $\ell:\bar B\times \bar B\to[0,\infty)$ denote the length metric induced by $F$; that is, $\ell(x,y)=\inf_c\int F(\dot c)$, where $c$ ranges over all smooth curves connecting $x$ with $y$.
Throughout, $(\bar B,F)$ is required to be a simple Finsler manifold.

For $v\in T\bar B$, let $\gamma_v:[t_-(v),t_+(v)]\to \bar B$ be the maximal geodesic with $\dot\gamma_v(0)=v$, so $\gamma_v(t_\pm(v))\in S^{n-1}$.
The geodesic flow on the unit tangent bundle $S\bar B:=\{v\in T\bar B:F(v)=1\}$ is thus given by
$$\Phi:\{(v,t)\in S\bar B\times\mathbb R:t_-(v)\le t\le t_+(v)\}\to S\bar B,\quad(v,t)\mapsto\phi^t(v)=\dot\gamma_v(t).$$
Moreover, let $\Gamma:=\{v\in S\bar B:\pi(v)\in S^{n-1},\,t_+(v)>0\}$ be the set of inward pointing unit vectors over the boundary, where $\pi:T\bar B\to \bar B$ denotes the footpoint projection.
Since $(\bar B,F)$ is simple, $t_+:\Gamma\to(0,\infty)$ is smooth, and% the restriction
$$\Phi:\{(v,t):v\in\Gamma,t\in(0,t_+(v))\}\to SB$$
is an orientation preserving diffeomorphism.

On $T\bar B\setminus\{0\}$, there is a natural one-form $\theta$, called Hilbert form:
$$T_wT\bar B\ni\xi\mapsto\theta_w(\xi)=g_w(w,D\pi(w)\xi)$$
It comes from the canonical one-form on $T^*\bar B$ via Legendre-transform; consequently, $d\theta$ is a symplectic two-form (cf.~\cite{Sh}, p.~26), and $\theta\wedge(d\theta)^{n-1}$ defines a volume form on $S\bar B$.
In fact, it is related to the Liouville form $\lambda$ via
$$\lambda=c_n\,\theta\wedge(d\theta)^{n-1},
\quad\textnormal{where}\quad
c_n:=\frac{(-1)^{n(n+1)/2+1}}{(n-1)!},$$
Hence, integration w.r.t.~Holmes-Thompson volume reads
$$\int_{\bar B}f\,d\mathrm{vol}=\frac{c_n}{\mathrm{vol}(S^{n-1})}\int_{S\bar B}f\circ\pi\;\theta\wedge(d\theta)^{n-1}\quad\forall\,f\in C(\bar B).$$
Both $d\theta$ and $\lambda$ are invariant w.r.t.~the geodesic flow (see \cite{Sh}, sect.~5.4).
\pagebreak[1]

Now, a Finsler version of Santal\'o's formula reads:

\begin{lem}\label{santalo1}
For every function $f\in L^1(S\bar B,\lambda)$, it holds
$$\int_{S\bar B}f\lambda=c_n\int_\Gamma\int_0^{t_+}f\circ\phi^t\,dt\,(d\theta)^{n-1}.$$
\end{lem}

\textsc{Proof:}
For $v\in S\bar B$, $t\in(t_-(v),t_+(v))$, fix some $(\xi,\tau)\in T_{(v,t)}(S\bar B\times\mathbb R)$.
Then it holds $\Phi^*\lambda_{(v,t)}
=c_n\cdot\big((\phi^t)^*\theta+dt\big)\wedge(d\theta)^{n-1}$, because of $(\phi^t)^*d\theta=d\theta$ and
\begin{eqnarray*}
\Phi^*\theta(\xi,\tau)
&=&\theta_{\phi^t(v)}(D\Phi(v,t)(\xi,\tau))
=\theta_{\phi^t(v)}(D\phi^t(v)\xi+\tau\tfrac d{dt}\phi^t(v))\\[1mm]
&=&(\phi^t)^*\theta(\xi)+\tau g_{\phi^t(v)}(\phi^t(v),\tfrac d{dt}\pi\circ\phi^t(v))
=((\phi^t)^*\theta+dt)(\xi,\tau).
\end{eqnarray*}
Hence, the claimed identity is obtained from transformation formula.
\hfill$\boxdot$\medskip

Further, set $S^{n-1}\times S^{n-1}\setminus\mathrm{diagonal}=:\Pi$ for shortness.
Then the map
$$\psi:\Gamma\to \Pi,\quad\psi(u)=(\pi(u),\exp(t_+(u)u))$$
is a diffeomorphism, w.r.t.~the orientation induced by $\psi$
This allows to express the integral in lemma~\ref{santalo1} in the following form:

\begin{prop}\label{santalo2}
The integral of any $f\in C(S\bar B)$ can be computed via
$$\int_{S\bar B}f\,d\lambda
=c_n\int_\Pi\int_0^\ell f\circ\phi^t\circ\psi^{-1}\,dt
\,(d_1d_2\ell)^{n-1}.$$
\end{prop}

\textsc{Proof:}
Set $V=\exp^{-1}(\bar B)\subset T\bar B$ and consider the map
$$\Psi:V\to \bar B\times \bar B,\quad\Psi(w)=(\pi(w),\exp(w)),$$
which is related to $\psi$ via $\psi(u)=\Psi(t_+(u)u)$, for $u\in\Gamma$.
Since all geodesics minimize distance, the first variation formula states that
$$d_1\ell(x,\exp(w))v=\frac{-g_w(w,v)}{F(w)}\quad\forall\,x\in \bar B,\;v,w\in T_x\bar B,\;w\in V\setminus\{0\}$$
Therefore,
\begin{equation}\label{Hilbertform}
\Psi^*(d_1\ell)_w=\pi^*d_1\ell(\pi(w),\exp(w))=\frac{-\theta_w}{F(w)}.
\end{equation}
Now, given $u\in\Gamma$, $\xi\in T_u\Gamma$, consider a smooth curve
$u:(-\varepsilon,\varepsilon)\to\Gamma$ with $u(0)=u$, $\dot u(0)=\xi$ and set $w(s)=t_+(u(s))u(s)$, thus $\psi\circ u=\Psi\circ w$.
From eqn.~\eqref{Hilbertform} and $g_u=g_w$ one infers
\begin{eqnarray*}
-\psi^*(d_1\ell)(\xi)\!
&=&\!-\Psi^*(d_1\ell)_w(\dot w(0))
=\frac1{F(w)}\theta_w(\dot w(0))\\
&=&\!\frac1{t_+(u)} g_w(t_+(u)u,\tfrac d{ds}\big|_0\pi\circ w(s))
=g_u(u,D\pi(u)\dot u(0))
=\theta(\xi).
\end{eqnarray*}
Using ``$d=d_1+d_2$'', one concludes that $d\theta=-d\psi^*d_1\ell=\psi^*d_1d_2\ell$.
Finally, prop.~\ref{santalo2} follows from lemma~\ref{santalo1} and the transformation formula.\hfill$\boxdot$\medskip\pagebreak[2]

To illustrate the geometric aspect of $D^2_{1,2}\ell$, consider $x,y\in\bar B$, $y\ne x$, and let $u=\frac{\exp_x^{-1}(y)}{\ell(x,y)}$ be the (Finsler) unit vector at $x$ pointing towards $y$ and
$$P_u:T_x\bar B\to T_x\bar B,\quad P_uv=v-g_u(u,v)u$$
the projection onto the $g_u$-orthogonal complement of $u$.
Then it holds:

\begin{prop}\label{D2elleq}
The mixed second derivative of $\ell(x,y)$ satisfies
$$D^2_{1,2}\ell(x,y)(v,w)
=\frac{-g_u(P_uv,D\exp_x^{-1}(y)w)}{\ell(x,y)}
\quad\forall\;v\in T_x\bar B,\;w\in T_y\bar B.$$
\end{prop}

\textsc{Proof:}
Let $c:(-\varepsilon,\varepsilon)\to\bar B\setminus\{x\}$ a smooth curve with $c(0)=y$ and $\dot c(0)=w$ and set $r(t)=\ell(x,c(t))$; hence one can write $c(t)=\exp_x(r(t)u(t))$ with $F(u(t))=1\;\forall\,t$.
Again, $d_1\ell(x,c(t))v=-g_{u(t)}(u(t),v)$ from the first variation formula, so
$$-D^2_{1,2}(x,y)(v,w)=\frac d{dt}\bigg|_0g_{u(t)}(u(t),v)=\frac{d^2}{2ds\,dt}\bigg|_0F(u(t)+sv)^2=g_u(v,\dot u(0)).$$
On the other hand,
$P_uD\exp_x^{-1}(y)w=P_u\big(\dot r(0)u(0)+r(0)\dot u(0)\big)=\ell(x,y)\dot u(0)$, because $g_u(u,\dot u(0))=\frac d{2dt}\big|_0F(u(t))^2=0$.\hfill$\boxdot$\medskip

{\small\textsc{Remarks.}
In case of $x,y\in S^{n-1}$ and $w\in T_yS^{n-1}$, one has $u(t)=\psi^{-1}(x,c(t))$, so $\dot u(0)=D_2\psi^{-1}(x,y)w$ and thus $-D^2_{1,2}\ell(x,y)(v,w)=g_{\psi^{-1}(x,y)}(v,D_2\psi^{-1}(x,y)w)$.
Further, if $(\xi_1,\ldots,\xi_{n-1},\upsilon_1,\ldots,\upsilon_{n-1})$ denote local coordinates on $\Pi$, the coordinate expression of $(d_1d_2\ell)^{n-1}$ reads
\begin{equation*}
(d_1d_2\ell)^{n-1}
=\Bigg(\sum_{i,j=1}^n\frac{\partial^2\ell}{\partial\xi_i\partial\upsilon_j}\cdot d\xi_i\wedge d\upsilon_j\Bigg)^{\!\!n-1}
=\det\bigg(\frac{\partial^2\ell}{\partial\xi_i\partial\upsilon_j}\bigg)\cdot(d\xi\wedge d\upsilon)^{n-1},
\end{equation*}
where $d\xi\wedge d\upsilon:=d\xi_1\wedge d\upsilon_1+\ldots+d\xi_{n-1}\wedge d\upsilon_{n-1}$.
Especially, the non-degeneracy of $d\theta$ implies that the determinant does not vanish.}

\section{An application to filling minimality.}

The Santal\'o-type integral formula from prop.~\ref{santalo2} can be used to obtain an equality between volume differences and certain integral of differences of boundary distances.
Again, set $\Pi=S^{n-1}\times S^{n-1}\setminus\mathrm{diagonal}$.

\begin{prop}\label{minfil}
Suppose $(\bar B,F)$ and $(\bar B,\tilde F)$ are simple Finsler manifolds with induced distances $\ell$ and $\tilde\ell$, respectively.
Then for the related Holmes-Thompson volumes, it holds
$$\mathrm{vol}_{\tilde F}(\bar B)-\mathrm{vol}_F(\bar B)=\frac{c_n}{\mathrm{vol}(S^{n-1})}\int_\Pi(\tilde\ell-\ell)\sum_{k=0}^{n-1}(dd_2\tilde\ell)^k\wedge(dd_2\ell)^{n-1-k}.$$
\end{prop}

\textsc{Proof:}
Taking $f\equiv1$ in prop.~\ref{santalo2}, one obtains
$$\mathrm{vol}_F(\bar B)=\frac{\mathrm{vol}_F(S\bar B)}{\mathrm{vol}(S^{n-1})}=\frac{c_n}{\mathrm{vol}(S^{n-1})}\int_\Pi\ell(d_1d_2\ell)^{n-1}.$$
Subtracting this from the corresponding expression for $\tilde F$ gives
$$\mathrm{vol}_{\tilde F}(\bar B)-\mathrm{vol}_F(\bar B)
=\frac{c_n}{\mathrm{vol}(S^{n-1})}\int_\Pi\tilde\ell\,(dd_2\tilde\ell)^{n-1}-\ell\,(dd_2\ell)^{n-1}.$$
The integrand can be decomposed into
$$\tilde\ell\,(dd_2\tilde\ell)^{n-1}-\ell\,(dd_2\ell)^{n-1}
=(\tilde\ell-\ell)\,(dd_2\tilde\ell)^{n-1}
+\ell\big((dd_2\tilde\ell)^{n-1}-(dd_2\ell)^{n-1}\big)$$
$$=(\tilde\ell-\ell)\,(dd_2\tilde\ell)^{n-1}
+\ell\,dd_2(\tilde\ell-\ell)\wedge\left(\sum_{k=0}^{n-2}(dd_2\tilde\ell)^k\wedge(dd_2\ell)^{n-2-k}\right).$$
Writing $\eta=\sum_{k=0}^{n-2}(dd_2\tilde\ell)^k\wedge(dd_2\ell)^{n-2-k}$ for simplicity, $\eta=1$ for $n=2$, while for $n>2$, $\eta$ is an exact $2(n-2)$-form of degree $n-2$ in each factor of $S^{n-1}\times S^{n-1}$.
Also, using ``$d=d_1+d_2$'' and ``$d_i^2=0$'', one obtains
$$\ell dd_2(\tilde\ell-\ell)-(\tilde\ell-\ell)dd_2\ell
=\ell dd_2\tilde\ell-\tilde\ell dd_2\ell
=d(\ell d_2\tilde\ell+\tilde\ell d_1\ell)-d_2\ell\wedge d_2\tilde\ell+d_1\ell\wedge d_1\tilde\ell.$$
Because $d_1\ell\wedge d_1\tilde\ell\wedge\eta$ and $d_2\ell\wedge d_2\tilde\ell\wedge\eta$ have degree $n$ in the first resp.~second variable, they cancel out.
For simplicity, set
$$(dd_2\tilde\ell)^{n-1}+dd_2\ell\wedge\eta=\sum_{k=0}^{n-1}(dd_2\tilde\ell)^k\wedge(dd_2\ell)^{n-1-k}=:\hat\eta,$$
so one infers from the above decomposition, that
$$\mathrm{vol}_{\tilde F}(\bar B)-\mathrm{vol}_F(\bar B)
=\frac{c_n}{\mathrm{vol}(S^{n-1})}\int_\Pi(\tilde\ell-\ell)\hat\eta+d(\ell d_2\tilde\ell\wedge\eta)+d(\tilde\ell d_1\ell\wedge\eta).$$
Herein, $\hat\eta$ is integrable, because $\hat\eta=n\cdot\int_0^1(dd_2((1-a)\ell+a\tilde\ell))^{n-1}\,da$ holds pointwise on $\Pi$, and the integrability of $(dd_2((1-a)\ell+a\tilde\ell))^{n-1}$ will be verified in cor.~\ref{d2ellintbil}.
\footnote{In fact, since two-forms can be muted without invoking sign changes, $\hat\eta$ can be considered a homogenous polynomial in $dd_2\ell$ and $dd_2\tilde\ell$ with all coefficients equal to 1.
The claimed integral representation thus follows from binomial expansion and the fact that
$$n\cdot{n-1\choose k}\,\int_0^1a^k(1-a)^{n-1-k}\,da=1\quad\forall\,k\in\{0,\ldots,n-1\}.$$}\\
Now, let $U(\varepsilon):=\{(x,y)\in \Pi:\ell(x,y)<\varepsilon\}$ denote a tubular $\varepsilon$-neighbour\-hood around $\partial \Pi=\mathrm{diag}(S^{n-1}\times S^{n-1})$.
Then Stokes' theorem implies
$$\int_{\Pi\setminus U(\varepsilon)}d(\ell d_2\tilde\ell\wedge\eta)
=\int_{\partial U(\varepsilon)}\ell d_2\tilde\ell\wedge\eta
=\varepsilon\int_{\partial U(\varepsilon)}d_2\tilde\ell\wedge\eta
=\varepsilon\int_{\Pi\setminus U(\varepsilon)}dd_2\tilde\ell\wedge\eta.$$
But $dd_2\tilde\ell\wedge\eta=\hat\eta-(dd_2\ell)^{n-1}$, so
$$\lim_{\varepsilon\to0}\int_{\Pi\setminus U(\varepsilon)}dd_2\tilde\ell\wedge\eta
=\int_\Pi\hat\eta-\int_\Pi(dd_2\ell)^{n-1}
\quad\textnormal{and}\quad\lim_{\varepsilon\to0}\varepsilon\!\int_{\Pi\setminus U(\varepsilon)}dd_2\tilde\ell\wedge\eta=0.$$
Likewise with $\tilde U(\varepsilon):=\{(x,y)\in \Pi:\tilde\ell(x,y)<\varepsilon\}$
$$\int_{\Pi\setminus\tilde U(\varepsilon)}d(\tilde \ell d_1\ell\wedge\eta)
=\varepsilon\int_{\Pi\setminus \tilde U(\varepsilon)}dd_1\ell\wedge\eta
=-\varepsilon\int_{\Pi\setminus \tilde U(\varepsilon)}dd_2\ell\wedge\eta\quad\stackrel{\varepsilon\to0}{\longrightarrow}\;0.$$
Thus, the integrals of the exact forms cancel, and one obtains the claimed equality.
\hfill$\boxdot$\medskip

\begin{cor}\label{minfilcor}
Let $(B,F)$ and $(B,\tilde F)$ be simple and such that $\ell(y,z)\le\tilde\ell(y,z)$ holds for all $y,z\in S^{n-1}$.
Then $\mathrm{vol}_F(\bar B)\le\mathrm{vol}_{\tilde F}(\bar B)$ with equality implying $\ell(y,z)=\tilde\ell(y,z)\;\forall\,y,z\in S^{n-1}$, provided one of the following criteria is satisfied:
\begin{enumerate}
\item
The dimension is $n=2$; or $n\le 4$ and there is a simple Finsler metric $\bar F$ having boundary distances $\bar\ell=\tilde\ell+\ell$.
\item
$\tilde\ell$ lies in an appropriate $C^2$-neighbourhood of $\ell$.
%\item
%There exists a $C^1$-homotopy $[0,1]\ni t\mapsto F_t$, with $F_0=F$ and $F_1=\tilde F$, of simple Finsler metrics, such that the corresponding boundary distances $\ell_t(x,y)$,  $x,y\in S^{n-1}$, are non-decreasing.
\end{enumerate}
\end{cor}

\textsc{Proof:}
Conditions 1 and 2 guarantee that $\hat\eta$ is a volume form, so prop.~\ref{minfil} implies the assertions:\\
For $n=2$, proposition \ref{santalo2} states that $\hat\eta=dd_2\ell+dd_2\tilde\ell$ corresponds to the sum of volume forms on the unit inward tangent bundle over $\partial B$, hence is again a volume form.
$\hat\eta$ can be decomposed as
$$\hat\eta=\tfrac12(dd_2\ell)^2
+\tfrac12\big(dd_2(\ell+\tilde\ell)\big)^2
+\tfrac12(dd_2\tilde\ell)^2
\quad\textnormal{for }n=3\textnormal{ and}$$
$$\hat\eta=\tfrac23(dd_2\ell)^3
+\tfrac13\big(dd_2(\ell+\tilde\ell)\big)^3
+\tfrac23(dd_2\tilde\ell)^3,\quad\textnormal{for }n=4.$$
Herein, the mixed term is a volume form, if $\ell+\tilde-\ell$ is simple.\pagebreak[1]

For the second condition: for every $\varepsilon>0$, there is a constant $\delta>0$, s.th.
$$|(dd_2\ell)^{n-1}(x,y)|\ge\delta|(dx\wedge dy)^{n-1}|\quad\textnormal{and}\quad\sup\|D^2_{1,2}\ell(x,y)\|\le\tfrac1\delta$$
for all $(x,y)\in \Pi$ with $\|x-y\|\ge\varepsilon$.
Accordingly, the expansion
\begin{eqnarray*}
\hat\eta&=&\sum_{k=0}^{n-1}(dd_2\ell)^{n-1-k}\wedge\big(dd_2\ell+dd_2(\tilde\ell-\ell))^k\\
&=&\sum_{i=0}^{n-1}\sum_{k=i}^{n-1}{k\choose i}(dd_2\ell)^{n-1-i}\wedge\big(dd_2(\tilde\ell-\ell))^i
\end{eqnarray*}
shows that $\hat\eta$ is dominated by the term $n(dd_2\ell)^{n-1}$, as long as $\|D^2_{1,2}(\tilde\ell-\ell)\|$ is smaller than some constant depending on $\delta$ and $n$.
Since $\varepsilon$ was arbitrary, this yields a $C^2$-neighbourhood for $\ell$ --- see however remark 2.
%If condition 3 is satisfied, one infers from the same expansion that $$\mathrm{vol}_{F_{t+s}}(B)-\mathrm{vol}_{F_t}(B)=\frac{c_n}{\mathrm{vol}(S^{n-1})}\int_\Pi(\ell_{t+s}-\ell_t)\big(n(dd_2\ell_t)^{n-1}+rem.\big)$$ is non-negative for $s>0$ sufficiently small, thus $\mathrm{vol}_{F_t}$ increases with $t$.
\hfill$\boxdot$\medskip

{\small\textsc{Remarks:}\\[1mm]
1. The first condition could be generalized for $n>4$.
Namely, one can choose $a_i\in[0,1]$ and $c_i>0$, such that $\hat\eta=\sum_ic_i\big(dd_2(a_i\ell+(1-a_i)\tilde\ell)\big)^{n-1}$.
Then $\hat\eta$ is a volume form, if $(1-a_i)\ell+a_i\tilde\ell$ are boundary distances of simple Finsler metrics $F_i$.\\[1mm]
2. Boundedness of $\|D^2_{1,2}(\tilde\ell-\ell)\|$ would require that $D^2_{1,2}\tilde\ell(x,y)$ and $D^2_{1,2}\ell(x,y)$ have the same asymptotic behaviour as $y\to x$; so $\tilde F$ and $F$ a priori would have to coincide on $S^{n-1}$ --- as was pointed out by S.~Ivanov.
Namely, due to prop.~\ref{D2elleq}, $D^2_{1,2}\ell$ becomes singular along the diagonal, indeed the scaling depends on direction.
To elude this deficiency, one can consider another criterion for positivity of $\hat\eta$ on $\{(x,y)\in\Pi:\|x-y\|<\varepsilon\}$, for $\varepsilon$ small.
Actually, in local coordinates $(\xi,\upsilon)$,
$$\hat\eta=n\!\int_0^1(dd_2((1-a)\ell+a\tilde\ell))^{n-1}\,da=n\!\int_0^1\det\bigg(\frac{\partial^2((1-a)\ell+a\tilde\ell)}{\partial\xi_i\partial\upsilon_j}\bigg)\,da\cdot(d\xi\wedge d\upsilon)^{n-1};$$
thus, it is sufficient to ensure that $\det\Big((1-a)\frac{\partial^2\ell}{\partial\xi_i\partial\upsilon_j}+a\frac{\partial^2\tilde\ell}{\partial\xi_i\partial\upsilon_j}\Big)$ does not vanish for $a\in(0,1)$.
In view of the remark after prop.~\ref{D2elleq}, this is satisfied, provided $\tilde\psi$ lies in a suitable $C^1$-neighbourhood of $\psi$ and $\tilde g_{\tilde\psi^{-1}(x,y)}$ is sufficiently close to $g_{\psi^{-1}(x,y)}$, for $(x,y)\in\Pi$.
In the remark after prop.~\ref{D2elldiag}, such a condition is stated in terms of $\ell,\tilde\ell$.}

\section{Analysis of $D^2_{1,2}\ell$ near the diagonal.}\label{diag.s}

Starting from prop.~\ref{D2elleq}, the objective of this section is to find two-sided estimates for $D^2_{1,2}\ell(x,y)$ as $x$ tends to $y$, in order to control the singularity of $(d_1d_2\ell)^{n-1}$ on the diagonal.

First, since $F$ is a Finsler metric, there is a constant $C_1>1$, such that
\begin{equation}\label{gijest}
\frac1{C_1^2}\|v\|^2\le g_u(v,v)\le C_1^2\|v\|^2\quad\forall\,u,v\in T_x\bar B,\;u\ne0
\end{equation}
where $\|v\|$ denotes the standard Euclidean norm on $\mathbb R^n$.
Furthermore, $C_1$ can be chosen independent of $x$, for compactness of $\bar B$.
As a consequence, one infers for the related distances
\begin{equation}\label{ellest}
\frac1{C_1}\|x-y\|\le\ell(x,y)\le C_1\|x-y\|\quad\forall\,x,y\in\bar B.
\end{equation}

The term $D\exp_x^{-1}(y)$ requires some scrutiny:
On a Finsler manifold, the exponential map at any point is known to be a local $C^1$-diffeomorphism on a neighbourhood of the origin, but of class $C^\infty$ only away from zero (see \cite{Sh}, thm.~11.1.1).
S.~Ivanov mentioned that the regularity is in fact $C^{1,1}$:

\begin{prop}\label{C11}
Let $(N,F)$ a smooth Finsler manifold without boundary.
Then for every point $p\in N$, the differential $D\exp_p$ of the exponential map is Lipschitz-continuous at $0\in T_pN$.
\end{prop}

\textsc{Proof (by S.~Ivanov):}
In local coordinates on a neighbourhood of $p$, the Finsler metric $F$ can be considered a function $F_1(x,v)$ of points $x\in\mathbb R^n$ and vectors $v\in\mathbb R^n$.
For simplicity, one can assume that $x(p)=0$ and extend $F_1$ arbitrarily to a smooth Finsler metric on the entire $\mathbb R^n$.
Define a family $F_t$, $t\in\mathbb R$ of ``blow-ups'' of the metric $F_1$ by $F_t(x,v) = F_1(tx,v)$.
This is a smooth family of metrics, so it defines a smooth family of exponential maps $E_t:\mathbb R^n\to\mathbb R^n$ (here $E_t$ is $\exp_0$ of the metric $F_t$).
More precisely, this family is smooth on any compact set separated away from the origin.
Let's consider it in a neighborhood of the unit sphere.
As $(\mathbb R^n,F_0)$ is a Minkowski space, $E_0$ is the identity, so its second derivative is zero.
Since $D^2E_t$ depends smoothly on $t$, there exists a constant $C>0$ such that, for $|t|$ small enough, $\|D^2 E_t(v)\| \le C|t|$ at any point $v$ of the unit sphere (here $\|\cdot\|$ is a norm on bilinear forms).
Because $F_1(tx,tv)=|t|\cdot F_t(x,v)$, the map $x\mapsto tx$ is a constant stretch and thus transfers geodesics in $(\mathbb R^n,F_t)$ to geodesics in $(\mathbb R^n,F_1)$.
Consequently, $E_1(v)=tE_t(v/t)$, so $DE_1(v)=DE_t(v/t)$ and $D^2E_1(v)=\frac1tD^2E_t(v/t)$ for all $t\ne 0$.
Rescaling back to the original metric, we get $\|D^2E_1(v)\|\le C$ for all $v$ on the sphere of radius $t>0$.
So $D^2E_1$ is bounded near the origin, hence $E_1=\exp_0$ is of class $C^{1,1}$.
\hfill$\boxdot$\medskip

Notice that, because of the smooth dependence of the generating vectorfield for the geodesic flow w.r.t.~changes in the Finsler metric, the corresponding maps $D^2E_{t,p}:S^{n-1}\to(\mathbb R^n\otimes\mathbb R^n)^*$ vary smoothly with $t$ and $p\in N$.
Therefore, the Lipschitz-constant $C$ can be chosen in a way that depends continuously on $p$.
This allows a uniform estimate in the next lemma:

\begin{lem}\label{D2ellest}
There exists a constant $C_2>1$, such that for all $x\not=y\in\bar B$ with $\|x-y\|<\frac1{C_2}$ and all $v\in T_x\bar B,w\in T_y\bar B$, it holds
$$\big|\ell(x,y)\cdot D^2_{1,2}\ell(x,y)(v,w)+g(P_uv,w)\big|
\;\le\;C_2\|x-y\|\cdot\|w\|\sqrt{g_u(P_uv,v)}.$$
\end{lem}

\textsc{Proof:}
When extending $F$ to a neighbourhood of $\bar B$, prop.~\ref{C11} guarantees the existence of some $C_3>0$, such that $\|D\exp_x(\tilde v)-\mathbf1\|\le C_3\|\tilde v\|$ for all $\tilde v\in T_x\bar B$ with $\|\tilde v\|<\frac1{C_3}$; and again $C_3$ can be selected independent of $x$, since $\bar B$ is compact.
If $\|\tilde v\|<\frac1{2C_3}$, then the inverse of $D\exp_x(\tilde v)$ satisfies
$$\|D\exp_x(\tilde v)^{-1}-\mathbf1\|\le\frac{\|D\exp_x(\tilde v)-\mathbf1\|}{1-\|D\exp_x(\tilde v)-\mathbf1\|}\le2C_3\|\tilde v\|$$
where the first inequality follows from
$$\|(A^{-1}-\mathbf1)w\|\le\|A-\mathbf1\|\cdot\|A^{-1}w\|\le\|A-\mathbf1\|\cdot\big(\|w\|+\|(A^{-1}-\mathbf1)w\|\big).$$
Taking $\tilde v=\exp_x^{-1}(y)$ in the above estimate, one obtains that
\begin{equation}\label{expest}
\|D\exp_x^{-1}(y)-\mathbf1\|\le2C_1^2C_3\cdot\|x-y\|
\end{equation}
as long as $\|x-y\|<\frac1{2C_1^2C_3}$, because due to ineqs.~\eqref{gijest} and~\eqref{ellest},
$$\|\exp_x^{-1}(y)\|\le C_1\cdot F(\exp_x^{-1}(y))=C_1\cdot\ell(x,y)\le C_1^2\|x-y\|.$$
Now, applying the Cauchy-inequality to the formula from prop.~\ref{D2elleq} states
$$\big|\ell(x,y)\cdot D^2_{1,2}\ell(x,y)(v,w)+g(P_uv,w)\big|
=\big|g_u(P_uv,(D\exp_x^{-1}(y)-\mathbf1)w)\big|$$
$$\le\sqrt{g_u((D\exp_x^{-1}(y)-\mathbf1)w,(D\exp_x^{-1}(y)-\mathbf1)w)}\cdot\sqrt{g_u(P_uv,v)}$$
for all $v\in T_x\bar B,w\in T_y\bar B$.
According to ineq.~\eqref{expest}, the first factor satisfies
\begin{eqnarray*}
\sqrt{g_u((D\exp_x^{-1}(y)-\mathbf1)w,(D\exp_x^{-1}(y)-\mathbf1)w)}
&\le&C_1\big\|(D\exp_x^{-1}(y)-\mathbf1)w\big\|\\
&\le&2C_1^3C_3\|x-y\|\cdot\|w\|
\end{eqnarray*}
provided that $\|x-y\|\le\frac1{2C_1^2C_3}$, which proves the assertion.
\hfill$\boxdot$\medskip

Restricting to the case of $x,y\in S^{n-1}$, let $e_{xy}\in T_xS^{n-1}$ denote the Euclidean unit vector tangent to the shortest arc on $S^{n-1}$ that connects $x$ with $y$.
Then $T_xS^{n-1}$ allows a decomposition into $\mathbb R\cdot e_{xy}$ and $T_{xy}:=T_xS^{n-1}\cap T_yS^{n-1}$, its orthogonal complement w.r.t.~the Euclidean scalar product $\langle\cdot\,,\cdot\rangle$.
The following estimates for $g_u(P_uv,w)$ will be needed in the sequel.

\begin{lem}\label{gijpest}
There exists a constant $C_4>C_2$, such that
\begin{eqnarray*}
g_u(P_ue_{xy},e_{xy})&\le&C_4\|x-y\|^2\;\ge\;g_u(P_ue_{yx},e_{yx}),\\
g_u(P_uv,v)&\ge&\frac1{C_4}\|v\|^2\qquad\forall\,v\in T_{xy}
\end{eqnarray*}
hold, as soon as $x,y\in S^{n-1}$ satisfy $0\ne\|x-y\|<\frac1{C_4}$.
\end{lem}

\textsc{Proof.}
First, when integrating ineq.~\eqref{expest} from the proof of lemma~\ref{D2ellest}, one obtains
\begin{eqnarray}
\|y-x-\exp_x^{-1}(y)\|
&=&\bigg\|\int_0^1\big(\mathbf1-D\exp_x^{-1}(ty+(1-t)x)\big)(y-x)\,dt\bigg\|\nonumber\\
&\le&\int_0^12C_1^2C_3\|ty-tx\|\cdot\|x-y\|\,dt=C_1^2C_3\|x-y\|^2.\label{xyest}
\end{eqnarray}
if $\|x-y\|\le\frac1{2C_1^2C_2}$.
On the other hand, one infers from plane geometry, that $$\bigg\|\frac{x-y}{\|x-y\|}-e_{yx}\bigg\|=\bigg\|\frac{y-x}{\|x-y\|}-e_{xy}\bigg\|=2\sin(s/4)\le\|x-y\|,$$
where $s=2\arcsin(\frac12\|x-y\|)$ is the Euclidean length of the shortest arc between $x$ and $y$ on $S^{n-1}$.
One concludes from the triangle inequality, that
$$\bigg\|\frac{\exp_x^{-1}(y)}{\|x-y\|}-e_{xy}\bigg\|\le(1+C_1^2C_3)\|x-y\|.$$
Since $P_u\Big(\frac{\exp_x^{-1}(y)}{\|x-y\|}\Big)=0$, one can apply ineqs.~\eqref{gijest} and~\eqref{xyest} to get
\begin{eqnarray*}
g_u(P_ue_{xy},e_{xy})
&=&g_u(P_u\bigg(e_{xy}-\frac{\exp_x^{-1}(y)}{\|x-y\|}\bigg),P_u\bigg(e_{xy}-\frac{\exp_x^{-1}(y)}{\|x-y\|}\bigg))\\
&\le&C_1^2\bigg\|\frac{\exp_x^{-1}(y)}{\|x-y\|}-e_{xy}\bigg\|^2
\le C_1^6C_3^2\|x-y\|^2.
\end{eqnarray*}
A similar reasoning would show the same estimate for $e_{yx}$, thereby verifying the first two claimed inequalities.

Next, let $z\in\mathbb R^n$ be the unique vector, s.th.~$g_u(v,z)=\langle v,y-x\rangle\;\forall\,v\in\mathbb R^n$.
Since $\langle v,y-x\rangle=0\;\forall\,v\in T_{xy}$, a Bessel-inequality reveals
$$1=g_u(u,u)\ge\frac{g_u(u,v)^2}{g_u(v,v)}+\frac{g_u(u,z)^2}{g_u(z,z)}\;\Rightarrow\;g_u(P_uv,v)\ge g_u(v,v)\frac{g_u(u,z)^2}{g_u(z,z)}.$$
For the numerator, the Cauchy-inequality and ineqs.~\eqref{xyest} and \eqref{ellest} show
\begin{eqnarray*}
g_u(u,z)
&=&\left\langle\frac{\exp_x^{-1}(y)}{\ell(x,y)},y-x\right\rangle
=\frac{\|y-x\|^2-\langle y-x-\exp_x^{-1}(y),y-x\rangle}{\ell(x,y)}\\
&\ge&\frac{\|x-y\|^2-\|x-y\|\cdot\|y-x-\exp_x^{-1}(y)\|}{\ell(x,y)}\\
&\ge&\frac{\|x-y\|^2(1-C_1^2C_2\|x-y\|)}{\ell(x,y)}
\ge\frac{\|x-y\|}{2C_1}.
\end{eqnarray*}
Further, ineq.~\eqref{gijest} implies a similar inequality for the dual metric $g_u^*$, so $g_u(z,z)=g_u^*((y-x)^T,(y-x)^T)\le C_1^2\|x-y\|^2$ in the denominator.
Collectively, these estimates demonstrate that $g_u(P_uv,v)\ge\frac1{4C_1^4}g_u(v,v)\ge\frac1{4C_1^6}\|v\|^2$ for all $v\in T_{xy}$.
At the end, $C_4$ can be chosen as the largest of the above constants.
\hfill$\boxdot$\medskip

Returning to the situation of prop.~\ref{minfil}, consider another simple Finsler metric $\tilde F$ on $\bar B$ with corresponding distance function $\tilde\ell$.

\begin{prop}\label{D2elldiag}
There exists a constant $C>1$, such that for arbitrary $a\in[0,1]$ and all $x\ne y\in S^{n-1}$ with $\|x-y\|\le\frac1{C}$, it holds:
$$\big|\big(d_1d_2\big((1-a)\ell+a\tilde\ell\,\big)(x,y)\big)^{n-1}\big|\le\frac{C}{\|x-y\|^{n-3}}\big|(dx\wedge dy)^{n-1}\big|.$$
\end{prop}

\textsc{Proof.}
Given $x,y\in S^{n-1}$, $-y\ne x\ne y$, let $e_1,\ldots,e_{n-2}$ be a basis of Euclidean unit vectors of $T_{xy}$, s.th.~$(e_1,\ldots,e_{n-2},e_{xy})$ and $(e_1,\ldots,e_{n-2},-e_{yx})$ form an oriented orthonormal basis of $T_xS^{n-1}$ and $T_yS^{n-1}$, respectively.
Then $d_1d_2\ell(x,y)^{n-1}=\det A\cdot(dx\wedge dy)^{n-1}$, where the coefficient matrix $A\in\mathbb R^{(n-1)\times(n-1)}$ has block shape
$$A=\left(\begin{array}{c|c}Q&c\\\hline r&s\end{array}\right)
\quad\textnormal{with}\quad\left\{\begin{array}{cccl}
q_{ij}&=&D^2_{1,2}\ell(x,y)(e_i,e_j)&(i,j\le n-2)\\
c_j&=&D^2_{1,2}\ell(x,y)(e_{xy},e_j)&(j\le n-2)\\
r_i&=&-D^2_{1,2}\ell(x,y)(e_i,e_{yx})&(i\le n-2)\\
s&=&-D^2_{1,2}\ell(x,y)(e_{xy},e_{yx})&
\end{array}\right.$$

Next, suppose that $\|x-y\|\le\frac1{C_2}$.
Then lemma~\ref{D2ellest} and ineqs.~\eqref{gijest}, \eqref{ellest} imply
$$\left|D^2_{1,2}\ell(x,y)(v,w)+\frac{g_u(P_uv,w)}{\ell(x,y)}\right|\le C_1^2C_2\|v\|\cdot\|w\|\quad\forall\,v\in T_x\bar B,w\in T_y\bar B.$$
Hence, the difference between the matrices $\frac{-1}{\ell(x,y)}(g_u(P_ue_i,e_j))_{i,j}$ and $Q$ is bounded by $C_1^2C_2$.
According to ineq.~\eqref{gijest}, the matrix $(g_u(P_ue_i,e_j))_{i,j}$ in turn is bounded from above by $C_1^2$.
Thus, due to ineq.~\eqref{ellest},
\begin{equation}\label{Best}
\|Q\|\le\frac{\|(g_u(P_ue_i,e_j))_{i,j}\|}{\ell(x,y)}+C_1^2C_2
\le\frac{C_1^3}{\|x-y\|}+C_1^2C_2.
\end{equation}
On the other hand, if $\|x-y\|\le\frac1{C_4}$, lemma~\ref{gijpest} states for all $v\in T_{xy}$, that $g_u(P_uv,v)$ is also bounded from below by $\frac{\|v\|^2}{C_4}$.
Hence, if $\|x-y\|\le\frac1{2C_1^3C_2C_4}$,
\begin{eqnarray}
v^TQv
&\ge&\frac{\|v\|^2}{C_4\ell(x,y)}-C_1^2C_2{\|v\|^2}
\ge\frac{\|v\|^2}{C_1C_4\|x-y\|}-C_1^2C_2{\|v\|^2}\nonumber\\
&\ge&\frac{\|v\|^2}{2C_1C_4\|x-y\|}\quad\forall\,v\in\mathbb R^{n-2}.
\label{Binvest}
\end{eqnarray}
Likewise, if $\|x-y\|\le\frac1{C_2}$, lemma~\ref{D2ellest} together with the Cauchy-inequality and ineq.~\eqref{ellest} show, that for $w\in T_yS^{n-1}$
\begin{eqnarray*}
\big|D^2_{1,2}\ell(x,y)(e_{xy},w)\big|
&\le&\frac{|g_u(P_ue_{xy},w)|+C_2\|x-y\|\cdot\|w\|\sqrt{g_u(P_ue_{xy},e_{xy})}}{\ell(x,y)}\\
&\le&C_1\sqrt{g_u(P_ue_{xy},e_{xy})}
\left(\frac{\sqrt{g_u(P_uw,w)}}{\|x-y\|}+C_2\|w\|\right).
\end{eqnarray*}
When $w=e_j$ and $\|x-y\|\le\frac1{C_4}$, one infers from lemma~\ref{gijpest} and ineq.~\eqref{gijest}:
\begin{equation}\label{exyw}
|c_j|\le C_1\sqrt{C_4}\|x-y\|\left(\frac{C_1\|e_j\|}{\|x-y\|}+C_2\right)\le C_1\sqrt{C_4}\left(C_1+\frac{C_2}{C_4}\right).
\end{equation}
Thus $\|c\|\le\sqrt{n-2}\,C_1\sqrt{C_4}\left(C_1+\frac{C_2}{C_4}\right)$, and the same estimate holds for $\|r\|$, too, since $\ell$ is symmetric when switching $x$ with $y$.
Also, setting $w=-e_{yx}$ and using lemma~\ref{gijpest} again, one obtains:
\begin{equation}\label{exyeyx}
|s|\le C_1\sqrt{C_4}\|x-y\|\left(\frac{\sqrt{C_4}\|x-y\|}{\|x-y\|}+C_2\right)=C_1(C_4+C_2\sqrt{C_4})\|x-y\|.
\end{equation}

After possibly taking larger constants, similar estimates like~\eqref{Best}--\eqref{exyeyx} hold true for the entries of $\tilde A$ corresponding to $\tilde\ell$, and even for the convex combination $\bar A:=(1-a)A+a\tilde A$ and its submatrices $\bar Q,\bar c,\bar r,\bar s$.
Especially, ineq.~\eqref{Binvest} states the claimed lower estimate for $D^2_{1,2}\bar\ell(x,y)$ on $T_{xy}$.
Now
$$\big(d_1d_2\big((1-a)\ell+a\tilde\ell\,\big)(x,y)\big)^{n-1}=\det\bar A\cdot(dx\wedge dy)^{n-1}.$$
As $\bar Q$ is invertible for $\|x-y\|$ sufficiently small, $\det\bar A$ can be computed via
\begin{equation}\label{detD2ell}
\det\bar A
=\det\left(\begin{array}{c|c}\bar Q&0\\\hline 0&1\end{array}\right)
\det\left(\begin{array}{c|c}\mathbf1&\bar Q^{-1}\bar c\\\hline \bar r&\bar s\end{array}\right)
=\det\bar Q\cdot(\bar s-\bar r\,\bar Q^{-1}\bar c),
\end{equation}
e.g.~by Laplace expansion in the last row.

Furthermore, one infers from ineqs.~\eqref{Best} and~\eqref{Binvest}, that
$$\det\bar Q\le\bigg(\!\frac{2C_1^3}{\|x-y\|}\!\bigg)^{\!\!n-2}\quad\textnormal{and}\quad
\|\bar Q^{-1}\|%\le\frac1{\min(\mathrm{Re}(\mathrm{spec}(\bar Q)))}
\le2C_1C_4\|x-y\|.$$
Combining the above estimates, eqn.~\eqref{detD2ell} implies for $\|x-y\|<\frac1{2C_1^3C_2C_4}$:
$$|\det\bar A|
\le|\det\bar Q|\cdot\big(|\bar s|+\|\bar r\|\cdot\|\bar c\|\cdot\|\bar Q^{-1}\|\big)\\
\le\frac{C}{\|x-y\|^{n-3}}$$
for some constant $C$, thereby proving the assertion.
\hfill$\boxdot$\medskip

{\small\textsc{Remarks.}\smallskip\\
1. The estimates in the proof also yield a sufficient condition for the non-vanishing of $\hat\eta=n\cdot\int_0^1(dd_2((1-a)\ell+a\tilde\ell))^{n-1}\,da$.
Namely, assume for $\|x-y\|<\varepsilon:=\frac1{2C_1^3C_2C_4}$, that
$$D^2_{1,2}(\tilde\ell-\ell)(x,y)(v,v)\le\frac{\varepsilon\|v\|^2}{2C_1C_4\|x-y\|}\quad\forall v\in T_{xy}$$
--- here $C_1,C_2,C_4$ are the constants related as before to $\ell$.
Then, in the above notations, $(1-a)Q+a\tilde Q$ is non-degenerate on $T_{xy}$, for all $a\in[0,1]$ and $\|x-y\|\le\varepsilon$.
Further, $(dd_2((1-a)\ell+a\tilde\ell))^{n-1}(x,y)=0$, if and only if
$$0=\frac{\det\bar A}{\det\bar Q}=s+a(\tilde s-s)-(r+a(\tilde r-r))(Q+a(\tilde Q-Q))^{-1}(c+a(\tilde c-c)).$$
Since $(d_1d_2\ell)^{n-1}$ is non-degenerate, $0\ne\det A$ and thus $0\ne s-rQ^{-1}c$.
Therefrom, one could deduce bounds on $|\tilde s-s|$, $\|\tilde r-r\|$, and $\|\tilde c-c\|$, that would guarantee $\det((1-a)A+a\tilde A)\ne0$ for all $a\in[0,1]$ and $\|x-y\|\le\varepsilon$.\smallskip\pagebreak[1]\\
2. In the model case of the Euclidean metric on $\bar B$, it follows from
$$D^2\ell(x,y)(v,w)=-\frac{\langle v,w\rangle-\langle v,u\rangle\cdot\langle u,w\rangle}{\|x-y\|},\quad u=\frac{y-x}{\|y-x\|},\quad e_{xy}=\frac{y-\langle x,y\rangle\,x}{\sqrt{1-\langle x,y\rangle^2}}$$
that $Q=-\|x-y\|^{-1}\cdot\mathbf1$, $r=c^T=(0,\ldots,0)$ and $s=\frac14\|x-y\|$.
This example might suggest, that $|s-rQ^{-1}c|\ge\frac1{C'}\|x-y\|$ should hold in general for some $C'>1$ and $\|x-y\|<\varepsilon$.
However, the estimates from lemma~\ref{D2ellest} and~\ref{gijpest} are too weak to verify this conjecture, since the error term is of the same order.}
\pagebreak[1]

The next corollary fills a gap in the proof of prop.~\ref{minfil}.

\begin{cor}\label{d2ellintbil}
$(dd_2((1-a)\ell+a\tilde\ell))^{n-1}$ is integrable on $\Pi$, $\forall\,a\in[0,1]$.
\end{cor}

\textsc{Proof.}
For continuity in the interior of $\Pi$, it is sufficient to verify integrablility in a neighbourhood of the diagonal.
To this end, let $z_k=y_k-x_k$; hence $(x_1,\ldots,x_n,z_1,\ldots,z_n)$ are new coordinates on $\mathbb R^n\times\mathbb R^n$, and the diagonal is just $\{(x,z):z=0\}$.
Further, $(dx\wedge dy)^{n-1}=(dx\wedge dz)^{n-1}$ plus a term that involves $(dx)^n$ and thus vanishes after restriction to $S^{n-1}\times S^{n-1}$.
Now $S^{n-1}\times S^{n-1}=\{(x,z):x\in S^{n-1},\,z\in S^{n-1}-x\}$, where $S^{n-1}-x$ is the sphere translated by $-x$.
One can switch from $z$ to polar-like coordinates $(r,\theta_1,\ldots,\theta_{n-2})$, with $r=\|z\|$ and local angle coordinates $(\theta_1,\ldots,\theta_{n-2})$ on $S_r^{n-1}\cap(S^{n-1}-x)$.
From transformation formula, there is a coefficient function $c=c(\theta)$ such that $(dz)^{n-1}=c(\theta)\cdot r^{n-2}dr\wedge(d\theta)^{n-2}$ on $S^{n-1}-x$.
Since $r=\|x-y\|$, one infers from prop.~\ref{D2elldiag} that
\begin{eqnarray*}
\big|\big(d_1d_2\big((1-a)\ell+a\tilde\ell\,\big)(x,y)\big)^{n-1}\big|
&\le&\frac{C}{\|x-y\|^{n-3}}\big|(dx\wedge dy)^{n-1}\big|\\
&=&C\cdot r\,\big|c(\theta)\,dr\wedge(d\theta)^{n-2}\wedge(dx)^{n-1}\big|
\end{eqnarray*}
holds for $(x,y)\in \Pi$ with $\|x-y\|<\frac1C$.
\hfill$\boxdot$

\section{A counterexample for positivity of $\hat\eta$.}

One could ask whether $\hat\eta$ (as defined in prop.~\ref{minfil}) is always a volume form in the given situation.
Unfortunately, this is wrong.

\begin{prop}\label{ctex}
There are simple Riemannian metrics, such that induced distances $\ell$ and $\tilde\ell$ satisfy $\tilde\ell(y,z)\ge\ell(y,z)\;\forall\,y,z\in S^{n-1}$, but s.th.~$\hat\eta$ is indefinite and there is no simple Finsler metric with boundary distances $\tilde\ell+\ell$.
\end{prop}

\textsc{Proof} by construction:\\
Let $\ell$ be the Euclidean distance on $\bar B\subset\mathbb R^3$.
Take $y=e_3=(0,0,1)$, $z=-e_3$, $v\in T_y\partial B$ and $w\in T_z\partial B$.
Using $v\perp e_3\perp w$, one obtains:
\begin{eqnarray*}
dd_2\ell(y,z)(v+0,0+w)&=&\frac{d^2}{ds\,dt}\bigg|_{s=t=0}\|y+sv-z-tw\|\\
&=&\frac d{dt}\bigg|_0\frac{\langle v,y-z-tw\rangle}{\|y-z-tw\|}=-\frac{\langle v,w\rangle}2.
\end{eqnarray*}
Further, let $\varphi:\bar B\to\bar B$ be a diffeomorphism with $\varphi(y)=y$, $\varphi(z)=z$, and consider the metric $\tilde\ell:=r\varphi^*\ell$ for some constant $r>1$.
Since $\tilde\ell$ is induced by the flat Riemannian metric $r^2\varphi^*\langle\cdot,\cdot\rangle$, $(\bar B,\tilde\ell)$ is still simple, and for $v,w\perp e_3$
\begin{eqnarray*}
&&dd_2\tilde\ell(y,z)(v+0,0+w)=\frac{d^2}{ds\,dt}\bigg|_{s=t=0}r\|\varphi(y+sv)-\varphi(z-tw)\|\\
&=&r\frac d{dt}\bigg|_0\frac{\langle D\varphi(y)v,y-\varphi(z-tw)\rangle}{\|y-\varphi(z-tw)\|}
=-r\frac{\langle D\varphi(y)v,D\varphi(z)w\rangle}2.
\end{eqnarray*}
Let $A,\mathbf1\in\mathbb R^{2\times 2}$ denote the matrices w.r.t.~$e_1,e_2$ of $D\varphi(y)^TD\varphi(z)$ and identity, resp.
The evaluation of $\hat\eta(y,z)$ on the basis of $T_{(y,z)}M\simeq T_y\partial B\oplus T_z\partial B$ given by $b_1=e_1+0,b_2=0+e_1,b_3=e_2+0,b_4=0+e_2$ reads
\begin{eqnarray*}&&\hat\eta(y,z)(b_1,b_2,b_3,b_4)\\[1mm]
&=&\tfrac12\big(dd_2\ell(y,z)^2+dd_2(\ell+\tilde\ell)(y,z)^2+dd_2\tilde\ell(y,z)^2\big)(b_1,b_2,b_3,b_4)\\[1mm]
&=&\tfrac12\big(\det(\tfrac12\mathbf1)+\det(\tfrac12\mathbf1+\tfrac r2A)+\det(\tfrac r2A)\big)\\
&=&\tfrac18\big(2+r\cdot\mathrm{tr}(A)+2r^2\det(A)\big).
\end{eqnarray*}

In order to get a negative result, $A$ should have two negative eigenvalues of different magnitude, so as to get a largely negative trace and a comparatively small but positive determinant.
A possible way to construct $\varphi$ with such kind of $A$ is to compose $\varphi$ of stretching the ball near $y,z$ with reciprocal factors and U-turn-torsion around the $e_3$-axis.

Therefore, consider the two parametrizations $$\psi_\pm:\mathbb R^2\to S^2\cap\{\pm x_3>0\},\quad\psi_\pm(\xi)=\frac{1}{\sqrt{1+\|\xi\|^2}}\left(\begin{array}{c}\xi_1\\\xi_2\\\pm1\end{array}\right)$$ for the upper and lower hemisphere.
Further, set $\rho(t)=\exp(-s^2t^2/2)$ for $s>1$ fixed and define maps $\phi_\pm:\mathbb R^2\to\mathbb R^2$ via $$\phi_\pm(\xi)=\left(\begin{array}{c}\xi_1\mp\rho(\xi_2)\xi_2/s\\\xi_2\pm s\rho(\xi_1)\xi_1\end{array}\right).$$
Finally, set $\varphi(x):=\|x\|\cdot\psi_\pm^{-1}\circ\phi_\pm\circ\psi_\pm\big(\frac x{\|x\|}\big)$ for $x_3\not=0$ and $\varphi(x)=x$ other\-wise.
Notice that $\varphi$ is differentiable along the equator, since $\psi_\pm^{-1}(x)=\frac{\pm1}{x_3}{x_1\choose x_2}$ and $\exp(-s^2x_{1,2}^2/2x_3^2)$ decays rapidly as $|x_3|\to0$.
The differential of $\phi_\pm$ is $$D\phi_\pm(\xi)=\left(\begin{array}{cc}1&\mp(1-s^2\xi_2^2)\rho(\xi_2)/s\\[1mm]\pm s(1-s^2\xi_1^2)\rho(\xi_1)&1\end{array}\right)$$ with $\det(D\phi_\pm(\xi))=1+(1-s^2\xi_1^2)\rho(\xi_1)(1-s^2\xi_2^2)\rho(\xi_2)$.
As follows from $\frac d{dt}(1-t)e^{-t/2}=\frac{t-3}2e^{-t/2}=0\Leftrightarrow t=3$, the coefficients $(1-s^2\xi_i^2)\rho(\xi_i)$ range between $-2e^{-3/2}$ and 1; so $\det(D\phi_\pm(\xi))\ge1-2e^{-3/2}>\frac12$.
Consequently $\phi_\pm$ are diffeomorphism, and thus $\varphi$ is also a diffeomorphism outside the origin, where it could be smoothened without loss of the boundary distance estimate.

Due to $D\psi_\pm(\pm e_3)=\mathbf1$, the matrix $A$ related to the specified $\varphi$ is $$A=D\phi_+(0)^TD\phi_-(0)=\left(\begin{array}{cc}1-s^2&s+1/s\\-s-1/s&1-1/s^2\end{array}\right)$$ and $\mathrm{tr}\,A=2-s^2-s^{-2}$, $\det A=4$.
For $\hat\eta(y,z)(b_1,b_2,b_3,b_4)$ be negative, it is then necessary that $$0>\tfrac18\big(2+r\cdot\mathrm{tr}(A)+2r^2\det(A)\big)=\tfrac18\big(2+r(2-s^2-s^{-2})+8r^2\big),$$
whereas $r$ must also fit to $s$ to guarantee that $r\varphi^*\ell>\ell$.
This in turn will hold, provided that $r\|D(\psi_\pm\circ\phi_\pm)(\xi)v\|\geq\|D\psi_\pm(\xi)v\|$ for all $\xi,v\in\mathbb R^2$.

Therefore, one computes
$$\|D\psi_\pm(\xi)v\|^2
=\left\|\frac d{dt}\bigg|_{t=0}\psi_\pm(\xi+tv)\right\|^2
=\frac{\|v\|^2}{1+\|\xi\|^2}-\frac{\langle v,\xi\rangle^2}{(1+\|\xi\|^2)^2},$$
$$\textnormal{so}\quad
\|D(\psi_\pm\circ\phi_\pm)(\xi)v\|^2
=\frac{\|D\phi_\pm(\xi)v\|^2}{1+\|\phi_\pm(\xi)\|^2}+\frac{\langle D\phi_\pm(\xi)v,\phi_\pm(\xi)\rangle^2}{(1+\|\phi_\pm(\xi)\|^2)^2}.$$
Applying $(a+b)^2\le2a^2+2b^2$ and the triangle inequality gives
\begin{eqnarray*}
\|\phi_\pm(\xi)\|^2&\le&\Bigg(\left\|\left(\begin{array}c\xi_1\\\pm s\xi_1\rho(\xi_1)\end{array}\right)\right\|+\left\|\left(\begin{array}c\mp\xi_2\rho(\xi_2)/s\\
\xi_2\end{array}\right)\right\|\Bigg)^2\\
&\le&2\xi_1^2(1+s^2\rho(\xi_1)^2)+2\xi_2^2(1+\rho(\xi_2)^2/s^2)\\
&\le&2\|\xi\|^2+2s^2\xi_1^2\rho(\xi_1)^2+2s^2\xi_2^2\rho(\xi_2)^2\\
&\le&2\|\xi\|^2+4
\end{eqnarray*}
because of $s^2\xi_i^2\rho(\xi_i)^2\le\frac{s^2\xi_i^2}{1+s^2\xi_i^2}<1$.
This states a bound for the quotient of the denominators:
$$\frac{1+\|\phi_\pm(\xi)\|^2}{1+\|\xi\|^2}\le\frac{5+2\|\xi\|^2}{1+\|\xi\|^2}\le5.$$

It remains to estimate the numerators.
In the sequel, vectors are interpreted as single-column-matrices, e.g.~$\langle v,w\rangle$ becomes $v^Tw$.
Then for $\xi\in\mathbb R^2$ fixed,
\begin{eqnarray*}
q(\xi)&:=&\sup_{v\in\mathbb R^2_*}\frac{\|v\|^2-(1+\|\xi\|^2)^{-1}\langle v,\xi\rangle}{\|D\phi_\pm(\xi)v\|^2-(1+\|\phi_\pm(\xi)\|^2)^{-1}\langle D\phi_\pm(\xi)v,\phi_\pm(\xi))^2}\\
&=&\sup_{v\in\mathbb R^2_*}\frac{v^T\big(\mathbf1-(1+\|\xi\|^2)^{-1}\xi\xi^T\big)v}{v^TD\phi_\pm(\xi)^T\big(\mathbf1-(1+\|\phi_\pm(\xi)\|^2)^{-1}\phi_\pm(\xi)\phi_\pm(\xi)^T\big)D\phi_\pm(\xi)v}.
\end{eqnarray*}
Since $\big(\mathbf1-(1+\|w\|^2)ww^T\big)^{-1}=\mathbf1+ww^T$ is positive and symmetric for all $w\in\mathbb R^3$, it has a unique positive, symmetric square root.
When substituting $v=D\phi_\pm(\xi)^{-1}\cdot
\sqrt{\mathbf1+\phi_\pm(\xi)\phi_\pm(\xi)^T}u$, one obtains
$$q(\xi)=\sup_{u\in\mathbb R^2_*}\frac{\big\|\sqrt{\mathbf1-(1+\|\xi\|^2)^{-1}\xi\xi^T}D\phi_\pm(\xi)^{-1}
\sqrt{\mathbf1+\phi_\pm(\xi)\phi_\pm(\xi)^T}u\big\|^2}{\|u\|^2}.$$
Writing $B(\xi)$ for the operator in the numerator, this is just the largest eigenvalue of $B(\xi)^TB(\xi)$.
It can be majorized by its trace; and using invariance of traces under cyclic permutation and linearity gives
\begin{eqnarray*}
&&q(\xi)<\mathrm{tr}(B(\xi)^TB(\xi))\\
&=&\mathrm{tr}\big(D\phi_\pm(\xi)^{-T}\big(\mathbf1-(1+\|\xi\|^2)^{-1}\xi\xi^T\big)D\phi_\pm(\xi)^{-1}
\big(\mathbf1+\phi_\pm(\xi)\phi_\pm(\xi)^T\big)\big)\\
&=&\mathrm{tr}\big(D\phi_\pm(\xi)^{-T}D\phi_\pm(\xi)^{-1}\big)
-(1+\|\xi\|^2)^{-1}\mathrm{tr}\big(D\phi_\pm(\xi)^{-T}\xi\xi^TD\phi_\pm(\xi)^{-1}\big)\\
&&+\;\mathrm{tr}\big(D\phi_\pm(\xi)^{-T}\big(\mathbf1-(1+\|\xi\|^2)^{-1}\xi\xi^T\big)
D\phi_\pm(\xi)^{-1}\phi_\pm(\xi)\phi_\pm(\xi)^T\big)\\
&\le&\mathrm{tr}\big(D\phi_\pm(\xi)^{-T}D\phi_\pm(\xi)^{-1}\big)\\
&&+\;\phi_\pm(\xi)^TD\phi_\pm(\xi)^{-T}\big(\mathbf1-(1+\|\xi\|^2)^{-1}\xi\xi^T\big)D\phi_\pm(\xi)^{-1}\phi_\pm(\xi)
\end{eqnarray*}
Because of $D\phi_\pm(\xi)^{-1}=\det(D\phi_\pm(\xi))^{-1}D\phi_\mp(\xi)$, the first summand reads
\begin{eqnarray*}
&&\mathrm{tr}\big(D\phi_\pm(\xi)^{-T}D\phi_\pm(\xi)^{-1}\big)
=\frac{\mathrm{tr}\big(D\phi_\mp(\xi)^TD\phi_\mp(\xi)\big)}{\det(D\phi_\pm(\xi))^2}\\
&=&\frac{2+s^2(1-s^2\xi_1^2)^2\rho(\xi_1)^2+s^{-2}(1-s^2\xi_2^2)^2\rho(\xi_2)^2}{\big(1+(1-s^2\xi_1^2)\rho(\xi_1)(1-s^2\xi_2^2)\rho(\xi_2)\big)^2}<4\big(s^2+3),
\end{eqnarray*}
due to $-\frac12<(1-s^2\xi_i^2)\rho(\xi_i)\le1$ as stated before.
Further, one can apply
$$\mathbf1-(1+\|\xi\|^2)^{-1}\xi\xi^T=\frac{\mathbf1+J\xi(J\xi)^T}{1+\|\xi\|^2},\quad\textnormal{with }J=\left(\begin{array}{cc}0&-1\\1&0\end{array}\right)$$
to rewrite the second summand and obtain
$$q(\xi)<4\bigg(s^2+3+\frac{\|D\phi_\mp(\xi)\phi_\pm(\xi)\|^2+\langle J\xi,D\phi_\mp(\xi)\phi_\pm(\xi)\rangle^2}{(1+\|\xi\|^2)}\bigg).$$
Now, $D\phi_\mp(\xi)\phi_\pm(\xi)
=\left(\begin{array}{c}\xi_1+\xi_1(1-s^2\xi_2^2)\rho(\xi_1)\rho(\xi_2)\\[1mm]\xi_2+\xi_2(1-s^2\xi_1^2)\rho(\xi_1)\rho(\xi_2)\end{array}\right)
+\left(\begin{array}{c}-s\xi_2^3\rho(\xi_2)^2\\[1mm]
s^3\xi_1^3\rho(\xi_1)^2\end{array}\right)$.\par
Because $\frac d{dt}t^m\rho(t)=(m-s^2t^2)t^{m-1}\rho(t)$ vanishes for ($t=0$ and) $t^2=ms^{-2}$, the functions $\xi_i^m\rho(\xi_i)$ have their maxima at $\big(\frac me\big)^{m/2}s^{-m}$.
Hence, the triangle inequality gives
$$\|D\phi_\mp(\xi)\phi_\pm(\xi)\|\le2\|\xi\|+\bigg(\frac3{2e}\bigg)^{3/2}\sqrt{1+s^{-2}}<2\|\xi\|+1,$$ and $(a+b)^2\le2a^2+2b^2$ implies $\|D\phi_\mp(\xi)\phi_\pm(\xi)\|^2\le8\|\xi\|^2+2$.
Also, 
\begin{eqnarray*}
\langle J\xi,D\phi_\mp(\xi)\phi_\pm(\xi)\rangle
&=&s^2(\xi_1\xi_2^3-\xi_1^3\xi_2)\rho(\xi_1)\rho(\xi_2)
+s\xi_2^4\rho(\xi_2)^2+s^3\xi_1^4\rho(\xi_1)^2\\
&\le&2\cdot3^{3/2}e^{-2}s^{-2}+2^2e^{-2}(s^{-3}+s^{-1})
<3.
\end{eqnarray*}
Assembling these estimates leads to
$$q(\xi)<4\bigg(s^2+3+\frac{8\|\xi\|^2+2+3^2}{(1+\|\xi\|^2)}\bigg)<4(s^2+14)$$
and shows that
$$\frac{\|D\psi_\pm(\xi)v\|}{\|D(\psi_\pm\circ\phi_\pm)(\xi)v\|}<10\sqrt{s^2+4}=:r\quad\forall\,v,\xi\in\mathbb R^2,\;v\ne0.$$
Finally, $s$ can be chosen sufficiently large to guarantee that
$$0>\hat\eta(y,z)(b_1,b_2,b_3,b_4)=\frac18\big(2+r(2-s^2-s^{-2})+8r^2\big).$$

This also proves that there must not be a simple Finsler metric with boundary distances $\ell+\tilde\ell$, because then $\hat\eta=\frac12(dd_2\ell)^2
+\frac12\big(dd_2(\ell+\tilde\ell)\big)^2
+\frac12(dd_2\tilde\ell)^2$ -- as a sum of volume forms --  would be positive.
\hfill$\boxdot$

\end{document}